
\baselineskip=14pt
\parskip=10pt
\def\halmos{\hbox{\vrule height0.15cm width0.01cm\vbox{\hrule height
  0.01cm width0.2cm \vskip0.15cm \hrule height 0.01cm width0.2cm}\vrule
  height0.15cm width 0.01cm}}

\magnification=\magstephalf

\def\1{{\overline{1}}}
\def\2{{\overline{2}}}
\parindent=0pt
\overfullrule=0in

\def\frac#1#2{{#1 \over #2}}
\centerline
{\bf
Automated Generation of Generating Functions Related to Generalized Stern's Diatomic Arrays 
 }
\centerline
{\bf
in the footsteps of Richard Stanley  
 }

\bigskip
\centerline
{\it Shalosh B. EKHAD and Doron ZEILBERGER}

\bigskip

{\bf Abstract}: Using {\it Symbolic Dynamic Programming} we describe algorithms, fully implemented in Maple,
for automatically generating generating functions introduced by Richard Stanley in his study
of generalized Stern arrays, generalized even further, to arrays defined in terms of general
sequences satisfying linear recurrences with constant coefficients, rather than just the Fibonacci
and k-bonacci sequences.

{\bf Appetizer: A Computational Challenge}

Mathematics is so sensitive to `{\it initial conditions}'. In other words the function

{\bf Mathematical Question} $\rightarrow$ {\bf Difficulty of Finding The Answer to that question}

is very `chaotic'. Just `tweaking' an easy question by a {\it proverbial} $\epsilon$ may change it from
{\it tractable} (and sometimes trivial) to  (very likely) {\it intractable} (but proving the
intractability for {\it sure} may also  be intractable!).

Here are a few of our favorite examples.

$\bullet$ We all know how to prove that there are infinitely many primes, but just
insert the word {\it twin} in front of  `primes' and no one (yet) can prove it.

$\bullet$ We all know how to compute the number of $10000$-step simple walks in the two dimensional square lattice ($4^{10000}$),
but stick {\it self-avoiding} in front of {\it walks} and we are all stumped.

$\bullet$ We all know, since Levi Ben Gerson, (exactly!) $700$ years ago, to compute the
number of {\it permutations} of size $10000$,  namely $10000!$, a certain 35660-digit integer
that Maple can compute (and display!) in $0.002$ seconds. Now stick the phrase
{\it $1234$-avoiding} in front of {\it permutations} and it would not take much longer,
(using the well-known second-order linear recurrence).
Now {\it transpose} the $2$ and the $3$ in ``$1234$-avoiding'', 
and we are all stumped on computing {\it exactly} the  integer that counts the number
of $1324$-avoiding permutations of length $10000$.

$\bullet$ The proof of the five-color theorem is half-page long, but the proof of the four-color theorem is
several-thousands-page long.

And we can come up with lots of other illustrations of the {\it sensitivity} of the ``difficulty'' function.

\vfill\eject

Here is yet another example.

{\bf Easy Problem}

Let 
$$
F_n(x) \,:=\, \prod_{i=0}^{n-1} \left ( 1 \,+\, x^{2^i}  \, + \, x^{2^{i+1}} \right ) \quad,
$$

and write
$$
F_n(x) \,  = \, \sum_{k \geq 0} \, a(n,k) \, x^k  \quad .
$$

Now let
$$
v(n) := \, \sum_{k \geq 0} \, a(n,k)^2 \quad .
$$

Find $v(10000)$.

Note that it is {\it hopeless} to compute this number from the definition. The degree, and hence the size, of the
polynomials $F_n(x)$ grow exponentially (in fact, it is $2(2^{n}-1)$).
Naively using the definition, we would have to add the squares of $2^{10001}-1$ numbers.

Using [S1], to be reviewed and extended later in this article, Maple can find $v(10000)$, a certain $6591$-digit number,
in a {\it split second.} 

Now comes our {\it tweak}

{\bf Hard Problem (at least for us)}

Let 
$$
G_n(x) \,:=\, \prod_{i=0}^{n-1} \left ( 1 \,+\, x^{2^i+1}  \, + \, x^{2^{i+1}+1} \right ) \quad,
$$

and write
$$
G_n(x) \,  = \, \sum_{k \geq 0} \, b(n,k) \, x^k  \quad  .
$$

Let
$$
w(n):=\sum_{k \geq 0} b(n,k)^2 \quad .
$$

Find the {\bf exact value} of $w(10000)$. 

One of us (DZ) is pledging a donation of $100$ US dollars to the OEIS (On-Line Encyclopedia of Integer Sequences)
in honor of the first (correct) solver of this very concrete problem.

Using the definition, Maple can compute the first $21$ values, starting at $n=0$:

$$
1, 3, 13, 55, 249, 1121, 5025, 22607, 101931, 460877, 2088687, 9482763, 43109307, 196163983, 
$$
$$
893222041, 4069162197, 18543631161, 84525140297, 385343891847, 1756959373157, 8011450183181 \quad .
$$

Unlike the previous problem, of computing $v(10000)$, for which it is easy to detect a simple ``pattern'' (see below) (and also easy to prove it
rigorously, also see below), and then easily deduce the $10000$-th term,
the modified problem, of computing $w(10000)$, seems much harder.

At any rate, the algorithms described later in the present article, that can handle everything in [S1] and [S2],
and much more,  fail miserably on this innocent problem. Of course, it is possible that there
{\bf exists} another algorithm that would make computing $w(10000)$ possible, but we have no clue,
and would love to know!

{\bf The original Stern Diatomic Array}

In the delightful article [S1], Richard Stanley's starting point was the double sequence $a(n,k)$ defined by
$$
F_n(x) \, = \, \sum_{k \geq 0} \, a(n,k) \, x^k \, = \, \prod_{i=0}^{n-1} \left ( 1 \,+\, x^{2^i}  \, + \, x^{2^{i+1}} \right ) \quad,
$$
and he was interested in the sequences (for positive integers $r$)
$$
u_r(n) \, := \, \sum_{k \geq 0} a(n,k)^r \quad .
$$
More generally for $\alpha=(\alpha_0, \dots, \alpha_{m-1})$ (where the $\alpha_i$ are non-negative integers), the
sequences
$$
u_\alpha(n) :=  \sum_{k \geq 0} a(n,k)^{\alpha_0}\, a(n,k+1)^{\alpha_1} \dots a(n,k+m-1)^{\alpha_{m-1}} \quad .
$$
He proved that all these sequences are $C$-finite, i.e., satisfy a {\it linear recurrence equation with constant coefficients},
or equivalently ([Z],[KP]) that their generating function 
$$
F_\alpha(x) \, := \, \sum_{n=0}^{\infty} u_\alpha(n) \, x^n \quad,
$$
is a {\bf rational function} of $x$.

In order to illustrate the general theory, he {\it humanly} proved that
$$
F_2(x)= \frac{1-2x}{1-5x+2x^2} \quad  .
$$

We will now redo, in {\it excruciating detail},  his proof, {\it our way}, in order to motivate the algorithm that will come later.
Our notation is  a little different than the one in [S1], but the bottom line is the same.

Since
$$
F_n(x) \, = \left ( 1 \,+\, x^{2^{n-1}}  \, + \, x^{2 \cdot 2^{n-1}} \right ) F_{n-1}(x) \quad,
$$
we have the recurrence that relates the entries in the $n$-row of Stern's array to those of the previous one.
$$
a(n,k) \, = \, a(n-1 \, , \, k) \,+ \,a(n-1 \, , \,  k- 2^{n-1}) \,+ \,  a(n-1 \, , \,  k - 2 \cdot 2^{n-1}) \quad .
$$

For reasons to be made clear later on, let's call $u_2(n)$, $f_{0}(n)$

Incorporating this in the definition of $u_2(n)$ (alias $f_{0}(n)$) we have
$$
f_{0} (n) \, = \, \sum_{k \geq 0}\, a(n,k)^2 \, = \, 
$$
$$
\sum_{k \geq 0}  
$$
$$
\, \left (a(n-1 \, , \, k) \, + \, a(n-1 \, ,\, k- 2^{n-1}) \, + \,a(n-1 \, , \, k- 2 \cdot 2^{n-1}) \right ) \cdot
$$
$$
\left (a(n-1 \, , \, k) \, + \, a(n-1 \, ,\, k- 2^{n-1}) \, + \,a(n-1 \, , \, k- 2 \cdot 2^{n-1}) \right ) 
$$
$$
\, = \, \sum_{k \geq 0}  \, a(n-1,k) \cdot a(n-1,k) 
\eqno(1.1)
$$
$$
+\, \sum_{k \geq 0}  \, a(n-1,k) \cdot a(n-1,k-2^{n-1})
\eqno(1.2)
$$
$$
+\, \sum_{k \geq 0}  \, a(n-1,k) \cdot a(n-1,k- 2 \cdot 2^{n-1})
\eqno(1.3)
$$
$$
+\, \sum_{k \geq 0}  \, a(n-1,k-2^{n-1}) \cdot a(n-1,k)
\eqno(1.4)
$$
$$
+\, \sum_{k \geq 0}  \, a(n-1,k-2^{n-1}) \cdot a(n-1,k-2^{n-1})
\eqno(1.5)
$$
$$
+\, \sum_{k \geq 0}  \, a(n-1,k-2^{n-1}) \cdot a(n-1,k- 2 \cdot 2^{n-1})
\eqno(1.6)
$$
$$
+\, \sum_{k \geq 0}  \, a(n-1,k- 2\cdot 2^{n-1}) \cdot a(n-1,k)
\eqno(1.7)
$$
$$
+\, \sum_{k \geq 0}  \, a(n-1,k-2 \cdot 2^{n-1}) \cdot a(n-1,k-2^{n-1})
\eqno(1.8)
$$
$$
+\, \sum_{k \geq 0}  \, a(n-1,k- 2 \cdot 2^{n-1}) \cdot a(n-1,k- 2 \cdot 2^{n-1}) \quad .
\eqno(1.9)
$$

Note that since the degree of $F_{n-1}(x)$ is $2(2^{n-1}-1)$, $(1.3)$  and $(1.7)$  are $0$.
Also note that by a `shift of the discrete variable $k$', $(1.1)$, $(1.5)$ and $(1.9)$ are the same.
By the commutativity of multiplication, $(1.2)$ and $(1.4)$ are identical, as are $(1.6)$ and $(1.8)$,
and by shifting the variable of summation, we see that these two pairs are identical to each other. Hence

$$
f_0(n) \, = \,3\, \sum_{k \geq 0}\, a(n-1 \, , \, k)^2  \, + \, \, 4\, \sum_{k \geq 0}\, a(n-1 \, , \, k) \, a(n-1 \, , \, k-2^{n-1})  \quad .
$$

The  first sum is $f_0(n-1)$ , but the second sum is a new creature, let's call if $f_1(n-1)$, where
$$
f_1(n) \, := \, \sum_{k \geq 0}\, a(n \, , \, k) \, a(n \,, \, k-2^{n})  \quad.
$$

So far we have
$$
f_0(n) \,= \, 3 f_0(n-1) + 4 f_1(n-1) \quad .
$$

We are {\bf forced} to consider $f_1(n)$.

We have

$$
f_{1} (n) \, = \, \sum_{k \geq 0}\, a(n \,, \, k)\, a(n \, ,  \, k - 2^n)\, = 
$$
$$
\, = \, \sum_{k \geq 0}  \,
\left( a(n-1 \,, \, k) + a(n-1, k-2^{n-1})+ a(n-1,k- 2 \cdot 2^{n-1}   \right ) \cdot
$$
$$
\left( a(n-1 \,, \, k- 2 \cdot 2^{n-1}) + a(n-1, k- 3 \cdot 2^{n-1})+ a(n-1,k- 4 \cdot 2^{n-1})   \right )
$$
$$
\, = \, \sum_{k \geq 0}  \, a(n-1,k) \cdot a(n-1,k- 2 \cdot 2^{n-1}) 
\eqno(2.1)
$$
$$
+\, \sum_{k \geq 0}  \, a(n-1,k) \cdot a(n-1,k- 3 \cdot 2^{n-1})
\eqno(2.2)
$$
$$
+\, \sum_{k \geq 0}  \, a(n-1,k) \cdot a(n-1,k- 4 \cdot 2^{n-1})
\eqno(2.3)
$$
$$
+ \, \sum_{k \geq 0}  \, a(n-1,k-2^{n-1}) \cdot a(n-1,k- 2 \cdot 2^{n-1}) 
\eqno(2.4)
$$
$$
+\, \sum_{k \geq 0}  \, a(n-1,k-2^{n-1}) \cdot a(n-1,k- 3 \cdot 2^{n-1})
\eqno(2.5)
$$
$$
+\, \sum_{k \geq 0}  \, a(n-1,k-2^{n-1}) \cdot a(n-1,k- 4 \cdot 2^{n-1})
\eqno(2.6)
$$
$$
+ \, \sum_{k \geq 0}  \, a(n-1,k -2 \cdot 2^{n-1} ) \cdot a(n-1,k- 2 \cdot 2^{n-1}) 
\eqno(2.7)
$$
$$
+\, \sum_{k \geq 0}  \, a(n-1,k - 2 \cdot 2^{n-1}) \cdot a(n-1,k- 3 \cdot 2^{n-1})
\eqno(2.8)
$$
$$
+\, \sum_{k \geq 0}  \, a(n-1,k - 2 \cdot 2^{n-1}) \cdot a(n-1,k- 4 \cdot 2^{n-1}) \quad .
\eqno(2.9)
$$

Once again since the  degree of $F_{n-1}(x)$ is $2(2^{n-1}-1)$, $(2.1)$, $(2.2)$, $(2.3)$, $(2.5)$, $(2.6)$ and $(2.9)$ vanish.
By shifting the summation variable $k$, both $(2.4)$ and $(2.8)$ equal $f_1(n-1)$ while $(2.7)$ is our old friend $f_0(n-1)$.

Hence we get, in addition to the previous equation, the following one:
$$
f_1(n) \,= \, f_0(n-1) + 2 f_1(n-1) \quad .
$$
Yea! We did not encounter any new `uninvited guests'. Defining  the
generating functions
$$
F_0(x)= \sum_{n=0}^{\infty} \, f_0(n) x^n \quad ,
$$
$$
F_1(x)= \sum_{n=0}^{\infty} \, f_1(n) x^n \quad,
$$

the above two recurrences translate to a {\bf system of two linear equations} in the two {\bf unknowns}
(also using  the {\it initial conditions} $f_0(0)=1$, $f_1(0)=0$)
$$
F_0(x)=1+ x(3 F_0(x)+ 4 F_1(x)) \quad, \quad  F_1(x)=0+ x( F_0(x)+ 2 F_1(x)) \quad.
$$
Of course, any  seventh-grader can solve this system, but why not use Maple?

Typing

{\tt latex(solve($\{$ F0=1+x*(3*F0+4*F1), F1=x*(F0+2*F1) $\}$,$\{$ F0,F1 $\}$ ));} 
we get
$$
 \left\{ {\it F0}=-{\frac {2\,x-1}{2\,{x}^{2}-5\,x+1}},{\it F1}={\frac {x}{2\,{x}^{2}-5\,x+1}} \right\}  \quad.
$$

Confirming that indeed $\sum_{n \geq 0} u_2(n) x^n$  equals $\frac{1-2x}{1-5x+2x^2}$ as claimed in [S1].
(Two lines below equation $(3)$ there).

Let's try and understand what is going on here. We started with an {\bf object of desire}, $f_0(n)$, and we were hoping to
relate it to $f_0(n-1)$. Alas, we were {\bf forced} to consider an {\it uninvited guest}, $f_1(n)$.
Analyzing $f_1(n)$, we were able to express it in terms of $f_1(n-1)$ and $f_0(n-1)$, and luckily, there
were no new `uninvited guests'. Taking the {\it $z$-transforms}, we got a system of two equations and two unknowns
and solving them, gave us the  generating function of $f_0(n)$, that we called $F_0(x)$, as well as
the generating function of $f_1(n)$, that we called $F_1(x)$. Of course we can
ungratefully disregard $F_1(x)$ at the end, if we wish, but we needed it in order to computer $F_0(x)$.

Let's now consider the general problem. Suppose that you have {\bf arbitrary} positive integers 
$$
0 \leq c_1 \leq c_2 \leq  \dots  \leq c_r \quad .
$$

We need the sequence
$$
\sum_{k \geq 0}  a(n,k+c_1) \, a(n,k+c_2) \, \cdots \, a(n,k+c_r) \quad .
$$
(Note that one can always take $c_1=0$.)

Using the recurrence 
$$
a(n,k)=a(n-1,k)+a(n-1,k-2^{n-1})+a(n-1,k-2 \cdot 2^{n-1}) \quad,
$$ 
above, this sum can be written as a sum of $3^r$ sums, each of them having the form
$$
f[d_1, \dots, d_r; \beta_1, \dots, \beta_r] (n-1) \, =
$$
$$
\sum_{k \geq 0}  a(n-1,k+d_1- \beta_1 2^{n-1}) \,\, a(n-1,k+d_2 - \beta_2 2^{n-1}) \,\, \cdots \, a(n-1,k+d_r- \beta_r 2^{n-1} ) \quad .
$$

Here $(d_1, \dots, d_r)$ is  a permutation of $(c_1, \dots, c_r)$ and $\beta_1, \dots, \beta_r$ are non-negative integers.

This {\bf forces} us to consider $f[d_1, \dots, d_r; \beta_1, \dots, \beta_r] (n)$ {\it in general}, not just the
initial case of $\beta_1=0, \dots, \beta_r=0$.

So we need to be able to express the general quantity
$$
f[d_1, \dots, d_r \,;\, \beta_1, \dots, \beta_r] (n):=
$$
$$
\sum_{k \geq 0}  a(n,k+d_1- \beta_1 2^{n}) \,\, a(n,k+d_2 - \beta_2 2^{n}) \, \, \cdots \,\, a(n,k+d_r- \beta_r 2^{n} ) \quad ,
$$
corresponding to the ``state'' $[d_1, \dots, d_r \,; \, \beta_1, \dots, \beta_r]$ in terms of other such ``states''.

We have to teach the computer:

$\bullet$ How to decide whether such a state is {\it dead on arrival}, i.e. identically zero (like  $(1.3)$, $(1.7)$, $(2.1$,
$(2.2)$, $(2.3)$, $(2.5)$,$(2.6)$, and $(2.9)$ in the example above).

$\bullet$ How to {\it automatically}  express each of these in terms of other such creatures.

Luckily, computer algebra comes to the rescue! Using the commutativity of multiplication, we can get a {\it canonical form} of each
state, by sorting the list of pairs 
$$
[d_1, \beta_1; d_2, \beta_2; \dots ; d_r, \beta_r] \quad,
$$
such that $\beta_1 \leq \beta_2 \leq \dots \leq \beta_r$ and $\beta_1=0$.
Of course as the $\beta_i$ change places, they must bring with them their corresponding $d_i$.

To each such state we  assign the {\bf monomial}
$$
x_1^{d_1} \cdots x_r^{d_r}  X_1^{\beta_1} \cdots X_r^{\beta_r} .
$$

First  one has to replace $X_i$ by $X_i^2$, since $k+d-\beta 2^{n-1}$ in the $n-1$ level becomes
$k+d-\beta 2^{n}\,=\, k+d- (2 \beta) 2^{n-1}$ at the $n$ level.

We let Maple multiply this (converted) monomial by
$$
\prod_{i=1}^{r} (1+X_i+X_i^2) \quad  ,
$$
and {\it expand} it, thereby expressing it as a sum of similar-looking monomials. 
We replace each monomial by its {\it canonical form} (sorting the powers of $X_i$ and permuting the corresponding
$x_i$ to move-along with their corresponding $X_i$, see Maple source-code).
Finally if the power of $X_1$ in the converted monomial is larger than $0$ we subtract
it from all the powers of $X_i$ (in other words if the power of $X_1$ is, $e$, we divide the monomial by
$X_1^{e} \cdots X_r^{e}$). This corresponds to shifting the variable $k$ in the sum corresponding to the given state.

Now each of these converted monomials corresponds to a state.

It is easy to see that there are only finitely many states (by the upper bound for the degree of $F_n(x)$),
so this process is {\it guaranteed} to {\it terminate}, and then
Maple automatically sets up a system of linear equations, that it can solve {\it all by itself}.

This approach works not just for the original Stern array, but for the more general scenario (introduced in [S1])
$$
F_n(x) \, = \, P(x) \prod_{i=0}^{n-1} Q( x^{b^i} ) \quad,
$$
for {\it any} polynomial $P(x)$ (before $P(x)$ was $1$), and {\it any} polynomial $Q(X)$, (before $Q(x)$ was $1+X+X^2$)
and for {\it any} integer $b \geq 2$ (before $b$ was $2$). 
The computer finds the equation for each
`still-to-do' state, by forming its corresponding monomial,
then replacing each $X_i$ by $X_i^{b}$ (due to the transition from the $n$ level to the $n-1$ level),
and then multiplying the resulting monomial by
$$
\prod_{i=1}^{r} Q(X_i)  \quad,
$$
expanding, taking the canonical forms, and discarding the monomials that
are `dead-on-arrival'.

This is fully implemented in procedure {\tt RS} in the Maple package {\tt StanleyStern.txt} available from the front of this article

{\tt https://sites.math.rutgers.edu/\~{}zeilberg/mamarim/mamarimhtml/stern.html} \quad ,

or directly from

{\tt https://sites.math.rutgers.edu/\~{}zeilberg/tokhniot/StanleyStern.txt}  \quad .

The function-call is

{\tt RS(P,Q,b,A,x,t)} \quad,

where {\tt P,Q,b} are as above and {\tt A} is $[\alpha_0, \dots, \alpha_{m-1}]$ featuring in Stanley's definition of $u_\alpha(n)$
mentioned above.

For example, to get the generating function of what Stanley [S1] calls $u_5(n)$, in the variable $t$
(rather than the variable $x$ that he is using), type

{\tt latex(RS(1+x+x**2,1,2,[5],x,t));} \quad \quad ,

immediately getting
$$
{\frac {20\,{t}^{2}+11\,t-1}{47\,{t}^{2}+14\,t-1}} \quad .
$$

Procedure {\tt RS} {\it dynamically} collects all the needed quantities  and by repeatedly invoking procedure {\tt Eq}
{\it dynamically} generates a system of linear equations, until it  does not encounter any {\it new guys}.

Of course, we know, {\it a priori} (as proved in [S1] and also follows from
our algorithm) that this must terminate, but even if we did not know that fact, we could have
set an {\bf upper limit} to the number of equations that we are willing to solve, and return {\tt FAIL} if
it got exceeded.

Since we have a {\it theoretical guarantee} that a rational function exists, and we can {\bf bound} the
degree of the denominator, we can use an {\it empirical} approach, and collect enough data
(as Stanley [S1] did in simple cases) and then fit it into the generating function (using something similar
to Maple's {\tt gfun[listtorec]}, but we prefer our own home-made version). This is implemented in procedure

{\tt RSe(p,q,b,A,N,x,t)}  \quad,

where $N$ is the number of data points used. Note that while for small cases, where the
order of the recurrences is expected to be relatively low,
this is much faster, but as noted in the above {\it appetizer}, this can't go very far,
since the polynomials $F_n(x)$ have exponential-size degree.

Nevertheless to compute generating functions for what Stanley calls $u_r(n)$ in [S1],
for small $r$ it works very fast. For example, to get the generating function of $u_{10}(n)$, type:

{\tt latex(RSe(1+x+x**2,1,2,[10],15,x,t));}

getting, in $0.16$ seconds
$$
-{\frac {4\,{t}^{4}+1852\,{t}^{3}+7945\,{t}^{2}+96\,t-1}{ \left( t+1 \right)  \left( 4\,{t}^{4}-200\,{t}^{3}-9601\,{t}^{2}-100\,t+1 \right) }} \quad .
$$

Using the non-guessing approach (using the algorithm outlined above), i.e. typing

{\tt latex(RS(1+x+x**2,1,2,[10],x,t));} \quad ,

gives that same thing, but in  $18.26$ seconds.

On the other hand, for the generating function of what is called $u_{11111}(n)$  in [S1], in other words, the generating function
of the sequence
$$
\sum_{k \geq 0} a(n,k)a(n,k+1)a(n,k+2)a(n,k+3)a(n,k+4) \quad,
$$

typing (using the {\it non-guessing} approach)

{\tt latex(RS(1+x+x**2,1,2,[1,1,1,1,1],x,t));}  \quad ,

yields, in $4.3$ seconds

$$
-4\,{\frac { \left( 4\,{t}^{4}-55\,{t}^{3}-69\,{t}^{2}-21\,t-3 \right) {t}^{2}}{ \left( t-1 \right) ^{3} \left( 47\,{t}^{2}+14\,t-1 \right) }} \quad ,
$$

while using the {\it guessing approach}, typing

{\tt latex(RSe(1+x+x**2,1,2,[1,1,1,1,1],20,x,t));}    \quad ,

yields the same thing in twice as long.

Eventually, the `guessing approach' will explode of course (as already mentioned before).

{\bf Sample Output files for the Maple package {\tt StanleyStern.txt}}

For the generating functions of [S1]'s $u_r(n)$ for $r \leq 25$, using the {\it guessing} (yet rigorous) approach see the
output file

{\tt https://sites.math.rutgers.edu/\~{}zeilberg/tokhniot/oStanleyStern1eA.txt}  \quad .

For the generating functions of [S1]'s $u_{1^r}(n)$, in other words of the sequences
$$
\sum_{k \geq 0} \, a(n,k) \, a(n,k+1) \, \cdots a(n,k+r-1) \quad,
$$

for $r \leq 9$, using the {\it guessing} (yet rigorous)  and non-guessing approaches respectively, see the
output files

{\tt https://sites.math.rutgers.edu/\~{}zeilberg/tokhniot/oStanleyStern2e.txt}  \quad ,

and

{\tt https://sites.math.rutgers.edu/\~{}zeilberg/tokhniot/oStanleyStern2.txt}  \quad  .

Note that the latter took quite a big longer.

{\bf Extension to Generalized Stern Arrays Induced by C-finite sequences}

In [S2], Richard Stanley extended his study to analogs where instead of $x^{b^i}$ that feature in the
definition of $F_n(x)$ above one has powers of the form $x^{F_i}$, and more generally, $x^{F^{(k)}_i}$
where $F_i$ are the Fibonacci numbers, and $F^{(k)}_i$ are the $k$-bonacci numbers defined by
$$
F^{(k)}_{i+1} \, = \, F^{(k)}_{i} \, + \, F^{(k)}_{i-1} \, + \, \dots \, +  F^{(k)}_{i-k+1} \quad,
$$
with initial conditions
$$
F^{(k)}_1 \, = \,F^{(k)}_2 \, =\, \dots \,= \,F^{(k)}_k \, = \,1 \quad .
$$

The sequences $\{b^i\}$, can be defined as a solution of the first-order recurrence
$$
f(i)\,=\, b \,f(i-1) \quad, \quad f(0)=1 \quad,
$$
and $F_i$, and $F^{(k)}_i$ are specific examples of $C$-finite sequences.

Let's first formally define a $C$-finite sequence.

{\bf Definition}: A $C$-finite sequence $f(i)$ (no relation to the Fibonacci numbers) of order $L$, is
a sequence defined by a recurrence of the form,
$$
f(i)=c_1 f(i-1) \, + \, c_2 f(i-2) \, + \, \dots \, + \, c_L f(i-L) \quad,
$$
where $c_1, \dots c_L$ are {\bf constants}, subject to some initial conditions
$$
f(0)=d_0 \, , \, f(1)=d_1  \, , \, \dots \, ,  \, f(L-1)= d_{L-1} \quad ,
$$
for some (other) constants $d_0, \dots d_{L-1}$.

We denote a $C$-finite sequence by the pair $[[d_0, \dots, d_{L-1}] \, ; \, [c_1, \dots, c_L]]$ and in this paper
assume that the $d's$ and $c's$ are integers, so the sequence $\{f(i)\}_{i=0}^{\infty}$ is an {\bf integer sequence}.

So in this notation, the sequence $\{2^i\}_{i=0}^{\infty}$ featuring in the original definition of the Stern array is 
denoted by $[[1],[2]]$ and more generally  $\{b^i\}_{i=0}^{\infty}$ is denoted by $[[1],[b]]$, while
the Fibonacci sequence is $[[0,1],[1,1]]$ and the $k$-bonacci sequence is $[[0,1^{k-1}], [1^k]]$.

Before going on, we need an important {\it observation}.

{\bf Observation}: Given a $C$-finite sequence $f(i)$ of order $L$, 
{\it any} linear combination of the form
$$
\sum_{i=M_1}^{M_2} a_i f(n+i) \quad,
$$
where $M_1<M_2$ are arbitrary integers, can be rewritten in the {\it canonical form}
$$
\sum_{i=0}^{L-1} a_i f(n+i) \quad ,
$$
where the new $a_i$ are of course different.

{\bf Proof}: $f(n+L)$ is a linear combination of $f(n), \dots, f(n+L-1)$. 
Let's prove, by induction on $J$, that $f(n+J)$ is a linear combination of $f(n), \dots, f(n+L-1)$, for all $J \geq L$.
It is true for $J=L$.
Assuming that $f(n+J)$ is a linear combination of $f(n), \dots, f(n+L-1)$, we get that
$f(n+J+1)$ is a linear combination of $f(n+1), \dots, f(n+L)$, and since $f(n+L)$ is a linear combination of $f(n), \dots, f(n+L-1)$,
it too is.  \halmos 

Now things are much more complicated, and we welcome the reader to study the source code of procedure {\tt RS} in the
other Maple package accompanying this article, {\tt SternCF.txt}, available from

{\tt https://sites.math.rutgers.edu/\~{}zeilberg/tokhniot/SternCF.txt}  \quad  .

The approach is similar, and let us describe it briefly.

We will consider the general problem of trying to find generating functions for the quantities $u_\alpha(n)$ but now
the array $a(n,k)$ is defined by an expression of the form
$$
F_n(x)= P(x) \prod_{i=0}^{n-1} \left ( c_0 + \sum_{j=1}^{M} c_j x^{ e_{[0,j]} f(i)+e_{[1,j]}f(i+1)+ \dots \, +\,e_{[L-1,j]}f(i+L-1)} \right ) \quad .
$$

Writing as before $F_n(x)=\sum_{k \geq 0}\,a(n,k)\,x^k$, 
we are interested in computing the generating function, if possible, of the quantity
$$
\sum_{k \geq 0}  a(n,k+d_1) \, a(n,k+d_2) \, \cdots \, a(n,k+d_r) \quad .
$$

Note that $a(n,k)$ satisfies the recurrence
$$
a(n,k)= c_0 a(n-1,k) \, + \, \sum_{j=1}^{M} c_j a \left (n-1, k- (e_{[0,j]} f(n-1)+e_{[1,j]}f(n)+ \dots \, +\,e_{[L-1,j]} f(n+L-2)) \right ) \quad .
$$

A typical `state' corresponds to the quantity
$$
\sum_{k \geq 0}  a(n,k+d_1- \beta_{1,0} f(n)- \beta_{1,1} f(n+1)- \dots - \beta_{1,L-1} f(n+L-1) ) \, \cdot
$$
$$
a(n,k+d_2- \beta_{2,0} f(n)- \beta_{2,1} f(n+1)- \dots - \beta_{2,L-1} f(n+L-1) ) \, \dots
$$
$$
a(n,k+d_r- \beta_{r,0} f(n)- \beta_{r,1} f(n+1)- \dots - \beta_{r,L-1} f(n+L-1) )  \quad .
$$

We denote this {\it state} by
$$
[[d_1; [\beta_{1,0},  \beta_{1,1}, \dots, \beta_{1,L-1}]] \quad,
$$
$$
[[d_2; [\beta_{2,0},  \beta_{2,1}, \dots, \beta_{1,L-1}]] \quad, 
$$
$$
\dots
$$
$$
[[d_r; [\beta_{r,0},  \beta_{r,1}, \dots, \beta_{r,L-1}]] \quad .
$$
Now, in addition to the variables $x_1, \dots x_r$ we have $r\,L$ variables
$$
X_{i,j} \quad , \quad 1 \leq i \leq r \quad, \quad 0 \leq j \leq L-1 \quad .
$$

The above state corresponds to the monomial
$$
x_1^{d_1} \cdots x_r^{d_r} \, \cdot \prod_{i=1}^r \prod_{j=0}^{L-1} X_{i,j}^{\beta_{i,j}} \quad .
$$

Before doing the `evolution' we must convert this monomial by translating from the $n$ level to the $n-1$ level,
(analogously to replacing $X_i$ by $X_i^{b}$ before). Now things are more complicated but the
computer does not mind (see procedure {\tt ROp}).

The {\it evolution equation} is obtained by multiplying this (adjusted) monomial by the polynomial
$$
\prod_{i=1}^{r} \left ( c_0 + \sum_{j=1}^{M} c_j X_{i,0}^{e_{[0,j]}} X_{i,1}^{e_{[1,j]}} \cdots X_{i,L-1}^{e_{[L-1 , j]}} \right ) \quad ,
$$

and expanding it, and then converting each of the monomials to its {\it canonical form}.

As before we build the system of equations dynamically, {\it except} that we are no longer guaranteed to terminate,
and indeed for many cases the process goes for ever. Hence we have another argument, {\tt LIMIT1}, telling us
to declare failure if the number of states (i.e. equations) exceeds it.

For the Fibonacci and $k$-bonacci cases considered in [S2], it always terminated, in the many cases we tried out,
as well as for many other cases. But not for {\it all} $C$-finite sequences!

For the $C$-finite sequence
$$
[[2,3],[3,-2]] \quad,
$$
alias $f(i)=\{ 2^i +1\}_{i=0}^{\infty}$ and
$$
F_n(x)= \prod_{i=0}^{n-1} \left ( 1 \,+\, x^{f(i)}  \, + \, x^{f(i+1)} \right ) \quad,
$$
even for {\tt A=[2]}, it seems  to never terminate.

For example entering

{\tt  RSmat([[1,1],[3,-2]],[[1,[0,0]], [1,[1,0]],[1,[0,1]]   ],1,x,[2],10000); }

shows that $10000$ does not suffice, and it is clear, that the set of states is {\it infinite}.
Hence the {\it computational challenge} at the  start of this article.

On the other hand for the Fibonacci and $k$-bonacci cases considered in [S2], it does seem to always terminate.

What makes them special?

Recall that a $PV$ number ({\it Pisot-Vijayaghavan} number) is a positive algebraic number that is larger than $1$ 
and such that all its conjugates  have absolute value less than $1$. We will call a $C$-finite sequence 
$$
[ [d_0, \dots d_{L-1}] \, , \, [c_1, \dots, c_L]]
$$

a $PV$-sequence if the largest root of the {\bf characteristic equation}
$$
X^L- c_1 x^{L-1} \, - \dots - c_{L-1}x-c_L \, = 0 \,  \quad,
$$
(the Golden ratio in the Fibonacci case) is a $PV$ number.

We believe that a careful study of our algorithm will be able to prove the following conjecture.

{\bf Conjecture}: Algorithm {\tt RS} of the Maple package {\tt SternCF.txt} (and its matrix version {\tt RSmat}) terminates
{\bf for all inputs} {\tt A}, (i.e. for computing the generating function of $u_\alpha(n)$), if and only if the $C$-finite sequence considered is a $PV$ sequence.

It is easy to see that not only the Fibonacci sequence, but the $k$-bonacci sequences are $PV$. That explains
why we were able to get answers for all the cases considered in [S2]. Of course, the system can get
very large, that's why we have a  matrix version of {\tt RS}, called {\tt RSmat}, that does not attempt to find the generating function
but only outputs the huge matrix $M$, and the vector $v$, such that the desired generating function is
$(I-tM)^{-1}v$, and it terminates in all the cases that we tried, as well as many other $PV$-sequences.

But the sequence $\{2^i+1\}$ that is `almost' $PV$, but not quite, seems to fail, hence this approach most probably
can {\bf not} prove that it is a rational generating function. This does not exclude the possibility that
another approach would prove that the sequence that we called $w(n)$, in the above challenge,
happens to be $C$-finite, but we strongly doubt it.
At any rate, if it  {\it is} $C$-finite its order must exceed $10$, while the order of the recurrence for what we called $v(n)$
above, is only two.

The Maple package {\tt SternCF.txt} produced quite a few output files, widely extending the computations in [S2].

For example, for generating functions for the sequences
$$
\sum_{k \geq 0} a(n,k)^r \quad,
$$
where 
$$
\sum_{k \geq 0} a(n,k) x^k = \prod_{i=1}^{n} \left ( 1+ x^{F_{i+1}} + x^{F_{i+2}} \right ) \quad ,
$$
for $1 \leq r \leq 6$ can be found in

{\tt https://sites.math.rutgers.edu/\~{}zeilberg/tokhniot/oSternCF1.txt}  \quad  .

The case $r=1$ is trivially $\frac{1}{1-3t}$, and the case $r=2$ confirms the generating function on the top of p.16 of [S2]
that Stanley probably got by the {\it guessing} approach. The degree of the denominator for the $r=3$ is already $35$
which means that for the guessing approach we would need to collect data up to $n=72$ which makes guessing impractical.

For the case $r=6$, the degree of the denominator (and numerator) of the generating function is $405$, that means that
we would need several `big bangs' to derive it by guessing.

Still with the same $a(n,k)$, but for the generating functions of
$$
\sum_{k \geq 0} a(n,k)a(n,k+1) \quad, \quad
\sum_{k \geq 0} a(n,k)a(n,k+1)a(n,k+2) \quad, \quad
\sum_{k \geq 0} a(n,k)a(n,k+1)a(n,k+2)a(n,k+3) \quad, \quad
$$
see the ouptput file

{\tt https://sites.math.rutgers.edu/\~{}zeilberg/tokhniot/oSternCF2.txt}  \quad  .

This took much longer, even though the degree was `only' $108$ for the last sequence
(but the number of states was much larger).

Moving right along to the Tribonacci sequence $T_i$ (alias $F^{(3)}_i$), and defining in analogy 

$$
\sum_{ k \geq 0} a(n,k) x^k = \prod_{i=1}^{n} (1+ x^{T_{i+1}} + x^{T_{i+2}} + x^{T_{i+3}} ) \quad ,
$$

the generating functions for $\sum_{k \geq 0} a(n,k)^r$ for $r=2$ and $r=3$ can be found here:

{\tt https://sites.math.rutgers.edu/\~{}zeilberg/tokhniot/oSternCF3.txt}  \quad  .

(the degree of the $r=3$ case is $567$).

The case $r=4$ is too big for us, but the `matrix version', {\tt RSmat}, that finds the matrix of coefficients
of the system, and enables computing many terms of the sequence, for  $\sum_{k \geq 0} a(n,k)^4$ can be found here:

{\tt https://sites.math.rutgers.edu/\~{}zeilberg/tokhniot/oSternCF3mat.txt}  \quad  .

The matrix in question has dimension $7245$.

Still with the same $a(n,k)$ (from the Tribonacci sequence), the generating function for $\sum_{k \geq 0} a(n,k)a(n,k+1)$ can be found here:

{\tt https://sites.math.rutgers.edu/\~{}zeilberg/tokhniot/oSternCF4.txt}  \quad  .

The generating function for $\sum_{k \geq 0} a(n,k)a(n,k+1)a(n,k+2)$ is too big (the system has $5004$ equations)
but we found the matrix that enabled us to compute the first $30$ terms, see here:

{\tt https://sites.math.rutgers.edu/\~{}zeilberg/tokhniot/oSternCF4mat.txt}  \quad  .

For the Quadonaci sequence, $Q_i$, (alias $F^{(4)}_i$), and defining in analogy
$$
\sum_{ k \geq 0} a(n,k) x^k = \prod_{i=1}^{n} \left ( 1+ x^{Q_{i+1}} + x^{Q_{i+2}} + x^{Q_{i+3}}  + x^{Q_{i+4}} \right ) \quad ,
$$
we only bothered to find the generating function of $\sum_{k \geq 0} a(n,k)^2$ that happens to have degree $504$. See here:

{\tt https://sites.math.rutgers.edu/\~{}zeilberg/tokhniot/oSternCF5.txt}  \quad  .

For the generating function of $\sum_{k \geq 0} a(n,k)a(n,k+1)$ that happens to have degree $1024$. See here:

{\tt https://sites.math.rutgers.edu/\~{}zeilberg/tokhniot/oSternCF6.txt}  \quad  .

For the analog of $\sum_{k \geq 0} a(n,k)^2$ for the $F^{(5)}_i$ case we  decided to only find
the $12751$-dimensional matrix, see

{\tt https://sites.math.rutgers.edu/\~{}zeilberg/tokhniot/oSternCF7mat.txt}  \quad  .

We also computed the quantities $J_r^{(k)}(t,x)$ considered  in section 5 of [S2] for quite a few $k$ and $r$.
(See [S2] for its definition.)

We believe that Conjecture 5.4 of [S2] is wrong as stated. Instead we have the

{\bf Corrected Conjecture 5.4 of [S2]}:
$$
J_3^{(k)}(t,x) \, = \,
$$
$$
\frac{
-{t}^{3\,k-3}{x}^{2\,k}{t}^{6}+ \left( {t}^{k-1} \right) ^{2}{x}^{k}{t}^{3}+{t}^{k-1}{x}^{k}{t}^{3}+2\,{t}^{3\,k-3}{x}^{2\,k}{t}^{3}+ \left( {t}^{k-1}
 \right) ^{2}{x}^{k}+{t}^{k-1}{x}^{k}-{t}^{3\,k-3}{x}^{2\,k}-1}{D_3^{(k)}(t,x)} \quad,
$$
where
$$
D_3^{(k)}(t,x) \, = \,
{t}^{3\,k-3}{x}^{2\,k+1}{t}^{9}- \left( {t}^{k-1} \right) ^{2}{x}^{k+1}{t}^{6}-{t}^{k-1}{x}^{k+1}{t}^{6}-{t}^{3\,k-3}{x}^{2\,k}{t}^{6}-{t}^{3\,k-3}{x}^{2\,k
+1}{t}^{6}
$$
$$
+ \left( {t}^{k-1} \right) ^{2}{x}^{k}{t}^{3}+ \left( {t}^{k-1} \right) ^{2}{x}^{k+1}{t}^{3}+{t}^{k-1}{x}^{k}{t}^{3}+{t}^{k-1}{x}^{k+1}{t}^{3}+2\,
{t}^{3\,k-3}{x}^{2\,k}{t}^{3}-{t}^{3\,k-3}{x}^{2\,k+1}{t}^{3}+{t}^{3}x
$$
$$
+ \left( {t}^{k-1} \right) ^{2}{x}^{k}- \left( {t}^{k-1} \right) ^{2}{x}^{k+1}+{t}^{k-
1}{x}^{k}-{t}^{k-1}{x}^{k+1}-{t}^{3\,k-3}{x}^{2\,k}+{t}^{3\,k-3}{x}^{2\,k+1}+x-1 \quad  ,
$$
and verified it for $k \leq 5$.

For the correct values of  $J_r^{(k)}(t,x)$ for $2 \leq r \leq 10$, and $2 \leq k \leq 4$ see the output file

{\tt https://sites.math.rutgers.edu/\~{}zeilberg/tokhniot/oSternCF9.txt}  \quad  .

For the conjectured values of  $J_r^{(k)}(1,x)$ for $2 \leq r \leq 20$, for {\it symbolic} (general) $k$, that match
the expressions given in [S2] for $r \leq 7$, see

{\tt https://sites.math.rutgers.edu/\~{}zeilberg/tokhniot/oSternCF10.txt}  \quad  .

Please note that these are still conjectures, but they were proved for $k \leq 6$ so they must be right.

Finally, Conjecture $5.6$ of [S2] is obviously wrong as stated, but if one replaces $a_i(t)$  by $a_i(t,t^{k-1})$ it
is probably possible to restate it correctly.

{\bf References}

[KP] Manuel Kauers  and Peter Paule, {\it ``The Concrete Tetrahedron''}, Springer, 2011.

[S1] Richard P. Stanley,
{\it Some linear recurrences motivated by Stern's Diatomic Array}, arXiv:1901.04647v1 [math.CO], 15 January 2019. 
{\tt https://arxiv.org/abs/1901.04647} \quad . \hfill\break
Also in: Amer. Math. Monthly {\bf 127} (2020), 99-111.

[S2]  Richard P. Stanley, {\it Theorems and conjectures on some rational generating functions},
arXiv:2101.02131v2 [math.CO], 11 January 2021. 
{\tt https://arxiv.org/abs/2101.02131} \quad . \hfill\break

[Z] Doron Zeilberger, {\it The C-finite Ansatz}, Ramanujan Journal {\bf 31} (2013), 23-32.\hfill\break
{\tt https://sites.math.rutgers.edu/\~{}zeilberg/mamarim/mamarimhtml/cfinite.html} \quad .

\vfill\eject

\bigskip
\hrule
\bigskip
Shalosh B. Ekhad and Doron Zeilberger, Department of Mathematics, Rutgers University (New Brunswick), Hill Center-Busch Campus, 110 Frelinghuysen
Rd., Piscataway, NJ 08854-8019, USA. \hfill\break
Email: {\tt [ShaloshBEkhad, DoronZeil] at gmail dot com}   \quad .

{\bf Exclusively published in the Personal Journal of Shalosh B. Ekhad and Doron Zeilberger and arxiv.org}

Originally Written: {\bf March 23, 2021}.

This version:  {\bf Oct. 20, 2024} [Correcting the values at the top 3, (pointed out by Guoce Xin, whom we thank)]

\end